\documentclass{article}
\usepackage{amssymb}

\input{tcilatex}

\begin{document}

\centerline{ Simplifying the axiomatization for the order affine geometry}

\centerline{Dafa Li} \centerline{Dept of mathematical sciences} \centerline
{Tsinghua University, Beijing 100084, China}

\centerline{Abstract}

Based on an ordering with directed lines and using constructions instead of
existential axioms, von Plato proposed a constructive axiomatization of
ordered affine geometry. There are 22 axioms for the ordered affine
geometry, of which the axiom I.7 is about the convergence of three lines
(ignoring their directions). In this paper, we indicate that the axiom I.7
includes much redundancy, and demonstrate that the complicated axiom I.7 can
be replaced with a simpler and more intuitive new axiom (called ODO) which
describes the properties of oppositely and equally directed lines. We also
investigate a possibility to replace the axiom I.6 with ODO.


\section{Introduction}

Heyting proposed the constructive axiomatization for the elementary geometry 
\cite{Heyting}. He adopted the concepts of distinct points and distinct
lines rather than the concepts of equal points and lines. von Plato proposed
the axioms of constructive geometry for which he introduced the concepts of
convergent lines instead of parallel lines and of apartness of a point from
a line instead of incidence of a point with a line \cite{jvp-95,Dalen}.
Recently, he presented a constructive theory of ordered affine geometry \cite%
{jvp-98}.

Automated theorem proving makes a great progress \cite{Wos, McCune-93,
McCune-96, McCune}. McCune proved that Robbins algebras are Boolean with the
theorem prover EQP \cite{McCune}. Since the problem was posed by Herbert
Robbins in the 1930s, it was conjectured that all Robbins algebras are
Boolean algebras.

By using theorem prover ANDP \cite{dli-97, dli-tableau}, we found a natural
deduction of Halting problem \cite{dli-halt}, and simplified von Plato's
axiomatization of constructive apartness geometry and orthogonality \cite%
{dli-00, dli-03}, and indicated that the equality axioms are not independent 
\cite{dli-04}.

In this paper, we simplify the axiomatization for ordered affine geometry by
using theorem prover. We indicate that the axiom 1.7 includes some
redundancy and show that the complicated axiom I.7 can be replaced with a
simpler and more intuitive new axiom.

\section{Simplifying the axiomatization for von Plato's order affine geometry%
}

The five basic relations DiPt, DiLn, Undir, L-Apt, and L-Con and four
constructions ln($a,b$), pt($a,b$), par$(l,m)$, and rev($l$) are used to
describe the axiomatization for the order affine geometry \cite{jvp-98}. In
total, von Plato's axiomatization has 22 axioms \cite{jvp-98}.

In this paper, we want to replace the axiom I.7 with a simpler and more
intuitive new axiom. Only four axioms I.5 to I.8 are concerned in this
paper. Only a basic relation Undir$(l,m)$ and only a construction rev$(l)$
appear in the axioms I.5 to I.8.

Undir$(l,m)$means $l$ and $m$ are unequally directed lines.

rev($l$) stands for the reverse of line $l$.

\subsection{List the four axioms I.5 to I.8 as follows.}

The axiom I.5, 
\begin{equation}
\thicksim Undir(l,l).
\end{equation}

The axiom I.6: 
\begin{equation}
Undir(l,m)\longrightarrow Undir(l,n)\vee Undir(m,n)  \label{I-6}
\end{equation}

The axiom I.7:%
\begin{eqnarray}
&&Undir(l,m)\&Undir(l,rev(m))  \nonumber \\
&\longrightarrow &Undir(l,n)\&Undir(l,rev(n))\vee
Undir(m,n)\&Undir(m,rev(n)).  \label{I-7}
\end{eqnarray}

The convergence $Con(l,m)$ is defined by means of the basic relation Undir\
as $Undir(l,m)\&Undir(l,rev(m))$ \cite{jvp-98}. Thus, $Con(l,m)$ means that $%
l$ and $m$ are convergent lines. Note that directions of lines $l$ and $m$
are ignored for $Con(l,m)$. And then by means of the definition of $Con(l,m)$%
, the axiom I.7 can be shortened as 
\begin{equation}
Con(l,m)\longrightarrow Con(l,n)\vee Con(m,n).  \label{conv}
\end{equation}

From Eq. (\ref{conv}), one can see that the axiom I.7 describes the
convergence of three lines (the direction is ignored)

The axiom I.8,

\begin{equation}
Undir(l,m)\vee Undir(l,rev(m))
\end{equation}

In \cite{jvp-98}, the new relation Inopp\ $(l,m)$ was defined as undir$%
(l,rev(m))$. By means of the definition of Inopp, the axiom I.8 can be
written as Undir$(l,m)$ $\vee $ Inopp$(l,m)$.

\subsection{The axiom I.7 is complicated and includes much redundancy}

The axiom I.7 can be rewritten equivalently as W1\&W2\&W3\&W4, where

\begin{eqnarray}
W1 &=&\thicksim Undir(l,m)\vee \thicksim Undir(l,rev(m))  \label{W-1} \\
&&\vee Undir(l,n)\vee Undir(m,n),  \nonumber \\
W2 &=&\thicksim Undir(l,m)\vee \thicksim Undir(l,rev(m))  \label{W-2} \\
&&\vee Undir(l,n)\vee Undir(m,rev(n)),  \nonumber \\
W3 &=&\thicksim Undir(l,m)\vee \thicksim Undir(l,rev(m))  \label{W-3} \\
&&\vee Undir(l,rev(n))\vee Undir(m,n),  \nonumber \\
W4 &=&\thicksim Undir(l,m)\vee \thicksim Undir(l,rev(m))  \label{W-4} \\
&&\vee Undir(l,rev(n))\vee Undir(m,rev(n)).  \nonumber
\end{eqnarray}

The classical concepts Opp and Dir are defined by means of the basic
relation Undir \cite{jvp-98}. $Opp(x,y)$ is defined as $\thicksim
Undir(x,rev(y))$. Thus, $Opp(x,y)$ means that the lines $x$ and $y$ are
oppositely directed lines. $Dir(x,y)$ is defined as $\thicksim Undir(x,y)$.
Thus, $Dir(x,y)$ means that $x$ and $y$ are equally directed lines.

By means of the definitions of Dir and Opp, W1, W2, W3, and W4 can be
rewritten as

\begin{eqnarray}
W1 &=&Dir(l,n)\&Dir(m,n)\longrightarrow Dir(l,m)\vee Opp(l,m),  \label{w1} \\
W2 &=&Dir(l,n)\&Opp(m,n)\longrightarrow Dir(l,m)\vee Opp(l,m),  \label{w2} \\
W3 &=&Opp(l,n)\&Dir(m,n)\longrightarrow Dir(l,m)\vee Opp(l,m),  \label{w3} \\
W4 &=&Opp(l,n)\&Opp(m,n)\longrightarrow Dir(l,m)\vee Opp(l,m).  \label{w4}
\end{eqnarray}

One can see that W1 in Eq. (\ref{w1}) means that if $l$ and $n$ are equally
directed line and $m$ and $n$ are equally directed lines, then lines $l$ and 
$m$ are equally directed lines or oppositely directed lines. One can know
that it is impossible for both $Dir(l,m)$ and $Opp(l,m)$ happen. By means of
the axiom I.6, one can know that lines $l$ and $m$ are equally directed
lines. Furthermore, from Eqs. (\ref{I-6}, \ref{W-1}), clearly the axiom I.6
implies W1. Thus, W1 is redundant.

W2 in Eq. (\ref{w2}) means that if $l$ and $n$ are equally directed line and 
$m$ and $n$ are oppositely directed lines, then lines $l$ and $m$ are
equally directed lines or oppositely directed lines. One can know that it is
impossible for both $Dir(l,m)$ and $Opp(l,m)$ happen. Factually, lines $l$
and $m$ should be oppositely directed lines. It means that W2 includes the
redundancy.

W3 in Eq. (\ref{w3}) means that if $l$ and $n$ are oppositely directed line
and $m$ and $n$ are equally directed lines, then lines $l$ and $m$ are
equally directed lines or oppositely directed lines. It is impossible for
both $Dir(l,m)$ and $Opp(l,m)$ happen. Factually, lines $l$ and $m$ should
be oppositely directed lines. It means that W3 includes the redundancy.

W4 in Eq. (\ref{w4}) means that if $l$ and $n$ are oppositely directed line
and $m$ and $n$ are oppositely directed lines, then lines $l$ and $m$ are
equally directed lines or oppositely directed lines. It is impossible for
both $Dir(l,m)$ and $Opp(l,m)$ happen. Factually, lines $l$ and $m$ should
be equally directed lines.

By other hand, the theorem prover ANDP finds a natural deduction of W4 from
only the axiom I-6. It means that the axiom I.6 implies W4 also. Thus, W4 is
redundant.

\subsection{The axiom I.7 can be replaced with the following simpler and
more intuitive new axiom called ODO.}

ODO:%
\begin{equation}
\thicksim Undir(l,rev(m))\&\thicksim Undir(l,n)\longrightarrow \thicksim
Undir(m,rev(n)).
\end{equation}

By means of the definitions of Opp and Dir, ODO becomes

\begin{equation}
Opp(l,m)\&Dir(l,n)\longrightarrow Opp(m,n),
\end{equation}%
which is just Theorem 3.5 \cite{jvp-98}. Thus, ODO means that if $l$ and $m$
are oppositely directed lines and $l$ and $n$ are equally directed lines,
then $m$ and $n$ are oppositely directed lines.

Theorem 1. Let the new axiomatization 1 be obtained from von Plato's one by
replacing axiom I.7 with ODO. Then, the new axiomatization 1 is equivalent
to von Plato's one.

Proof. It is known that ODO can be derived from von Plato's axiomatization.
Conversely, we can derive the axiom I.7 from the new axiomatization 1 as
follows.

To reduce the difficulty to derive the axiom I.7, we consider W1, W2, W3,
and W4 in Eqs. (\ref{W-1}, \ref{W-2}, \ref{W-3}, \ref{W-4}). In stead of
finding a deduction of the axiom I.7, we only need to find deductions of W1,
W2, W3, and W4, respectively. After lots of trials, the natural deductions
of W1, W2, W3, and W4 are obtained as follows.

W1 is derived from only the axiom I.6 without using other axioms. Ref.
Appendix A.

W2 is derived from the axioms I.5, I.6, and ODO without using other axioms.
Ref. Appendix B.

W3 is derived from the axioms I.5, I.6, ODO without using other axioms. Ref.
Appendix C.

W4 is derived from only the axiom I.6 without using other axioms. Ref.
Appendix D.

\section{Automated natural deduction}

The the natural deduction is adapted from Gentzen system. We push
quantifiers inside as possible as we can in our theorem prover ANDP. Thus,
we can have more chances to apply rules for propositional logic. ANDP uses
the following rules \cite{dli-97}.

The rules for propositional logic:

MP (modus ponens), MT (modus tollens), IMP ($\thicksim A\vee B$\ from $%
A\rightarrow B$), LDS ($B$ from $A$ and $\thicksim A\vee B$), RDS ($B$ from $%
A$ and $B\vee \thicksim A$),\ CP (conditional proof), SIMPlication ($A$ or $%
B $ from $A\wedge B$), CASES (or dilemma).

The rules for quantifiers:

US (universal specialization), UG (universal generalization), EG
(existential quantifier), EE (eliminate existential quantifier), SUB
(substitute a term $t$ for an individual variable $v$).

The wffs of the axioms I.5, I.6, I.7, I.8, and ODO are listed as follows.

The axiom I.5, 
\begin{equation}
(\forall x)\thicksim Undir(x,x)
\end{equation}

The axiom I.6

\begin{equation}
(\forall x)(\forall y)[Undir(x,y)\rightarrow (\forall z)[Undir(x,z)\vee
Undir(y,z)]]
\end{equation}

For the axiom I.7, the wffs of W1, W2, W3, and W4 are given as follows. 
\begin{eqnarray}
W1 &=&(\forall x)(\forall y)(\forall z)[\thicksim Undir(x,y)\vee \thicksim
Undir(x,rev(y)) \\
&&\vee Undir(x,z)\vee Undir(y,z)]  \nonumber \\
W2 &=&(\forall x)(\forall y)(\forall z)[\thicksim Undir(x,y)\vee \thicksim
Undir(x,rev(y)) \\
&&\vee Undir(x,z)\vee Undir(y,rev(z))]  \nonumber \\
W3 &=&(\forall x)(\forall y)(\forall z)[\thicksim Undir(x,y)\vee \thicksim
Undir(x,rev(y)) \\
&&\vee Undir(x,rev(z))\vee Undir(y,z)]  \nonumber \\
W4 &=&(\forall x)(\forall y)(\forall z)[\thicksim Undir(x,y)\vee \thicksim
Undir(x,rev(y)) \\
&&\vee Undir(x,rev(z))\vee Undir(y,rev(z))]  \nonumber
\end{eqnarray}

The axiom I.8,

\begin{equation}
(\forall x)(\forall y)[Undir(x,y)\vee Undir(x,rev(y))]
\end{equation}

ODO:%
\begin{equation}
(\forall x)(\forall y)(\forall z)[\thicksim Undir(x,rev(y)]\wedge \thicksim
Undir(x,z)\rightarrow \thicksim Undir(y,rev(z))]
\end{equation}

Note that only one predicate (i.e. Undir$(l,m)$) and only one function
symbol (i.e. rev$(l)$) appear in the axioms I.5 to I.8, and ODO. For our
automated prover, Undir$(l,m)$ is written as Undir $l$ $m$, and the function
symbol rev$(l)$ is written as [rev $l$]. The quantifiers $(\forall x)$ and $%
(\exists x)$ are written as $(Ax)$ and $(Ex)$, respectively. The connectives 
$\rightarrow $, $\vee $, and $\wedge $ are written as $->$, $|$, and $\&$,
respectively. \textquotedblleft $\symbol{126}$\textquotedblright\ stands for
\textquotedblleft $\thicksim $\textquotedblright\ (i.e. \textquotedblleft
not\textquotedblright ).

\section{Replacing the axiom I.6 with ODO}

It is also interesting to replace the axiom I.6 with ODO.

Theorem 2. Let the axiomatization 2 be obtained from von Plato's
axiomatization by replacing axiom I.6 with ODO. Then, the axiomatization 2
is equivalent to von Plato's axiomatization.

Proof. It is known that ODO can be derived from von Plato's axiomatization.
Thing left to do is to derive the axiom I.6 from the axiomatization 2. Our
theorem proving ANDP finds a natural deduction of the axiom I.6 from the
axioms I.7, I.8, and ODO. Ref. Appendix E.

\section{Appendix A. The axiom I.6 implies W1}

Solution

1. (Ax)(Ay)[UNDIR x y -$>$ (Az)[UNDIR x z $|$ UNDIR y z]] ASSUMED-PREMISE

2. UNDIR v1 v2 \ \ \ \ \ \ \ \ \ \ \ \ \ \ \ \ \ \ \ \ \ \ \ \ \ \ \ \ \ \ \
\ \ \ \ \ \ \ \ \ ASSUMED-PREMISE

3. \symbol{126}UNDIR v1 v3 \ \ \ \ \ \ \ \ \ \ \ \ \ \ \ \ \ \ \ \ \ \ \ \ \
\ \ \ \ \ \ \ \ \ \ \ \ ASSUMED-PREMISE

4. (Ay)[UNDIR v1 y -$>$ (Az)[UNDIR v1 z $|$ UNDIR y z]] \ \ \ \ US (v1 x) 1

5. UNDIR v1 v2 -$>$ (Az)[UNDIR v1 z $|$ UNDIR v2 z] \ \ \ \ \ \ \ US (v2 y) 4

6. (Az)[UNDIR v1 z $|$ UNDIR v2 z] \ \ \ \ \ \ \ \ \ \ \ \ \ \ \ \ \ \ \ \ \
\ \ \ \ \ \ \ \ \ \ \ \ MP 5 2

7. UNDIR v1 v3 $|$ UNDIR v2 v3 \ \ \ \ \ \ \ \ \ \ \ \ \ \ \ \ \ \ \ \ \ \ \
\ \ \ \ \ \ \ \ \ \ \ \ \ \ \ US (v3 z) 6

8. UNDIR v2 v3 \ \ \ \ \ \ \ \ \ \ \ \ \ \ \ \ \ \ \ \ \ \ \ \ \ \ \ \ \ \ \
\ \ \ \ \ \ \ \ \ \ \ \ \ \ \ \ \ \ \ \ \ \ \ \ \ \ \ \ \ \ \ LDS 7 3

9. UNDIR v2 v3 \ \ \ \ \ \ \ \ \ \ \ \ \ \ \ \ \ \ \ \ \ \ \ \ \ \ \ \ \ \ \
\ \ \ \ \ \ \ \ \ \ \ \ \ \ \ \ \ \ \ \ \ \ \ \ \ \ \ \ \ \ \ SAME 8

10. \symbol{126}UNDIR v1 v3 -$>$ UNDIR v2 v3 \ \ \ \ \ \ \ \ \ \ \ \ \ \ \ \
\ \ \ \ \ \ \ \ \ \ \ \ \ \ \ \ CP 9

11. UNDIR v1 v3 $|$ UNDIR v2 v3 \ \ \ \ \ \ \ \ \ \ \ \ \ \ \ \ \ \ \ \ \ \
\ \ \ \ \ \ \ \ \ \ \ \ \ \ \ IMP 10

12. UNDIR v1 [rev v2]-$>$ UNDIR v1 v3 $|$ UNDIR v2 v3 \ \ \ \ \ \ \ \ \ \ CP
11

13. [\symbol{126}UNDIR v1 [rev v2] $|$ UNDIR v1 v3] $|$ UNDIR v2 v3 \ \ \ \
\ \ \ \ IMP 12

14. UNDIR v1 v2

-$>$ [\symbol{126}UNDIR v1 [rev v2] $|$ UNDIR v1 v3] $|$ UNDIR v2 v3 \ \ \ \
\ \ \ \ \ CP 13

15. [[\symbol{126}UNDIR v1 v2 $|$ \symbol{126}UNDIR v1 [rev v2]] $|$ UNDIR
v1 v3]

$|$ UNDIR v2 v3 \ \ \ \ \ \ \ \ \ \ \ \ \ \ \ \ \ \ \ \ \ \ \ \ \ \ \ \ \ \
\ \ \ \ \ \ \ \ \ \ \ \ \ \ \ \ \ \ \ \ \ \ \ \ \ \ \ \ \ \ \ \ \ \ \ \ IMP
14

16. (Ax)(Ay)(Az)[[[\symbol{126}UNDIR x y $|$ \symbol{126}UNDIR x [rev y]]

$|$ UNDIR x z] $|$ UNDIR y z] \ \ \ \ \ \ \ \ \ \ \ \ \ \ \ \ \ \ \ \ \ \ \
\ \ \ \ \ \ \ \ \ \ \ \ \ \ \ \ \ \ \ \ \ \ \ \ \ \ UG 15

17. (Ax)(Ay)[UNDIR x y -$>$ (Az)[UNDIR x z $|$ UNDIR y z]]

-$>$ (Ax)(Ay)(Az)[[[\symbol{126}UNDIR x y $|$ \symbol{126}UNDIR x [rev y]]

$|$ UNDIR x z] $|$ UNDIR y z] \ \ \ \ \ \ \ \ \ \ \ \ \ \ \ \ \ \ \ \ \ \ \
\ \ \ \ \ \ \ \ \ \ \ \ \ \ \ \ \ \ \ \ \ \ \ \ \ \ \ CP 16

\section{Appendix B. Derive W2 from the axioms I.5, I.6, and ODO}

To reduce the difficulty to derive W2, first we derive the following result
from the axiom I.5 and ODO. Ref. Appendix B.1. 
\begin{equation}
(\forall x)(\forall y)[Undir(x,rev(y))\rightarrow Undir(y,rev(x))]
\label{oo}
\end{equation}

Then, W2 can be derived from I.5, I.6, and the above result. Ref. Appendix
B.2.

\subsection{B.1. Derive Eq. (\protect\ref{oo}) from the axiom I.5 and ODO}

Show that

(Ax)(Ay)[UNDIR x [rev y] -$>$ UNDIR y [rev x]]

from the following premises:

(Ax)(Ay)(Az)[\symbol{126}UNDIR x [rev y] \& \symbol{126}UNDIR x z -$>$ 
\symbol{126}UNDIR y [rev z]];

\symbol{126}(Ex)UNDIR x x.

Solution

1. (Ax)(Ay)(Az)[\symbol{126}UNDIR x [rev y] \& \symbol{126}UNDIR x z

-$>$ \symbol{126}UNDIR y [rev z]] \ \ \ \ \ \ \ \ \ \ \ \ \ \ \ \ \ \ \ \ \
\ \ \ \ \ \ \ \ \ \ \ \ \ \ \ \ \ \ \ \ \ \ PREMISE

2. \symbol{126}(Ex)UNDIR x x \ \ \ \ \ \ \ \ \ \ \ \ \ \ \ \ \ \ \ \ \ \ \ \
\ \ \ \ \ \ \ \ \ \ \ \ \ \ \ \ \ \ \ \ \ \ PREMISE

3. UNDIR v1 [rev v2] \ \ \ \ \ \ \ \ \ \ \ \ \ \ \ \ \ \ \ \ \ \ \ \ \ \ \ \
\ \ \ \ \ \ \ ASSUMED-PREMISE

4. \symbol{126}UNDIR v3 v3 \ \ \ \ \ \ \ \ \ \ \ \ \ \ \ \ \ \ \ \ \ \ \ \ \
\ \ \ \ \ \ \ \ \ \ \ \ \ \ \ \ \ \ \ \ \ \ \ US (v3 x) 2

5. (Ay)(Az)[\symbol{126}UNDIR v4 [rev y] \& \symbol{126}UNDIR v4 z-$>$ 
\symbol{126}UNDIR y [rev z]]

\ \ \ \ \ \ \ \ \ \ \ \ \ \ \ \ \ \ \ \ \ \ \ \ \ \ \ \ \ \ \ \ \ \ \ \ \ \
\ \ \ \ \ \ \ \ \ \ \ \ \ \ \ \ \ \ \ \ \ \ \ \ \ \ \ \ \ \ \ \ \ \ \ \ \ US
(v4 x) 1

6. (Az)[\symbol{126}UNDIR v4 [rev v5] \& \symbol{126}UNDIR v4 z

-$>$ \symbol{126}UNDIR v5 [rev z]] \ \ \ \ \ \ \ \ \ \ \ \ \ \ \ \ \ \ \ \ \
\ \ \ \ \ \ \ \ \ \ \ \ \ \ \ \ \ \ \ \ US (v5 y) 5

7. \symbol{126}UNDIR v4 [rev v5] \& \symbol{126}UNDIR v4 v6

-$>$ \symbol{126}UNDIR v5 [rev v6] \ \ \ \ \ \ \ \ \ \ \ \ \ \ \ \ \ \ \ \ \
\ \ \ \ \ \ \ \ \ \ \ \ \ \ \ \ \ \ \ US (v6 z) 6

8. \symbol{126}UNDIR v4 [rev v1] \& \symbol{126}UNDIR v4 v2

-$>$ \symbol{126}UNDIR v1 [rev v2] \ \ \ \ \ \ \ \ \ \ \ \ \ \ \ \ \ \ \ \ \
\ \ \ \ \ \ \ \ \ \ \ \ \ \ \ \ \ \ \ \ SUB 7

9. \symbol{126}[\symbol{126}UNDIR v4 [rev v1] \& \symbol{126}UNDIR v4 v2] \
\ \ \ \ \ \ \ \ \ \ MT 3 8

10. UNDIR v4 [rev v1] $|$ UNDIR v4 v2 \ \ \ \ \ \ \ \ \ \ \ \ \ \ \ \ \ \ \
\ \ \ \ \ DE.MORGAN 9

11. \symbol{126}UNDIR v2 v2 \ \ \ \ \ \ \ \ \ \ \ \ \ \ \ \ \ \ \ \ \ \ \ \
\ \ \ \ \ \ \ \ \ \ \ \ \ \ \ \ \ \ \ \ \ \ \ \ \ \ \ \ \ \ SUB 4

12. UNDIR v2 [rev v1] $|$ UNDIR v2 v2 \ \ \ \ \ \ \ \ \ \ \ \ \ \ \ \ \ \ \
\ \ \ \ \ \ SUB 10

13. UNDIR v2 [rev v1] \ \ \ \ \ \ \ \ \ \ \ \ \ \ \ \ \ \ \ \ \ \ \ \ \ \ \
\ \ \ \ \ \ \ \ \ \ \ \ \ \ \ \ \ \ \ \ \ \ RDS 11 12

14. UNDIR v2 [rev v1] \ \ \ \ \ \ \ \ \ \ \ \ \ \ \ \ \ \ \ \ \ \ \ \ \ \ \
\ \ \ \ \ \ \ \ \ \ \ \ \ \ \ \ \ \ \ \ \ \ SAME 13

15. UNDIR v1 [rev v2] -$>$ UNDIR v2 [rev v1] \ \ \ \ \ \ \ \ \ \ \ \ \ \ \ \
CP 14

16. (Ax)(Ay)[UNDIR x [rev y] -$>$ UNDIR y [rev x]] \ \ \ \ \ \ \ \ \ \ UG 15

\subsection{B.2. Derive W2 from the axioms I.5, I.6, and Eq. (\protect\ref%
{oo}):}

Show that

[(Ax)\symbol{126}UNDIR x x

\& (Ax)(Ay)[UNDIR x [rev y]-$>$ UNDIR y [rev x]]]

\& (Ax)(Ay)[UNDIR x y -$>$ (Az)[UNDIR x z $|$ UNDIR y z]]

-$>$ (Ax)(Ay)(Az)[[[\symbol{126}UNDIR x y $|$ \symbol{126}UNDIR x [rev y]]

$|$ UNDIR x z] $|$ UNDIR y [rev z]]

Solution

1. [(Ax)\symbol{126}UNDIR x x

\& (Ax)(Ay)[UNDIR x [rev y] -$>$ UNDIR y [rev x]]]

\& (Ax)(Ay)[UNDIR x y -$>$ (Az)[UNDIR x z $|$ UNDIR y z]] \ \ 

\ \ \ \ \ \ \ \ \ \ \ \ \ \ \ \ \ \ \ \ \ \ \ \ \ \ \ \ \ \ \ \ \ \ \ \ \ \
\ \ \ \ \ \ \ \ \ \ \ \ \ \ \ \ \ \ \ \ \ \ \ \ \ \ \ \ \ \ \ \ \ \ \ \ \ \
\ \ \ \ ASSUMED-PREMISE

2. \symbol{126}(Ex)UNDIR x x \ \ \ \ \ \ \ \ \ \ \ \ \ \ \ \ \ \ \ \ \ \ \ \
\ \ \ \ \ \ \ \ \ \ \ \ \ \ \ \ \ \ \ \ \ \ \ \ \ \ \ \ \ \ \ \ \ \ \ SIMP 1

3. (Ax)(Ay)[UNDIR x [rev y] -$>$ UNDIR y [rev x]] \ \ \ \ \ \ \ \ \ \ \ \ \
\ \ SIMP 1

4. (Ax)(Ay)[UNDIR x y -$>$ (Az)[UNDIR x z $|$ UNDIR y z]] \ \ \ SIMP 1

5. UNDIR v1 [rev v2] \ \ \ \ \ \ \ \ \ \ \ \ \ \ \ \ \ \ \ \ \ \ \ \ \ \ \ \
\ \ \ \ \ \ \ \ \ \ \ \ \ \ \ \ \ \ \ \ \ \ \ \ \ \ \ ASSUMED-PREMISE

6. \symbol{126}UNDIR v1 v3 \ \ \ \ \ \ \ \ \ \ \ \ \ \ \ \ \ \ \ \ \ \ \ \ \
\ \ \ \ \ \ \ \ \ \ \ \ \ \ \ \ \ \ \ \ \ \ \ \ \ \ \ \ \ \ \ \ \ \ \ \
ASSUMED-PREMISE

7. (Ay)[UNDIR v1 y -$>$ (Az)[UNDIR v1 z $|$ UNDIR y z]] \ \ \ \ \ \ \ \ \ US
(v1 x) 4

8. UNDIR v1 [rev v2]

-$>$ (Az)[UNDIR v1 z $|$ UNDIR [rev v2] z] \ \ \ \ \ \ \ \ \ \ \ \ \ \ \ \ \
\ \ US (rev(v2) y) 7

9. (Az)[UNDIR v1 z $|$ UNDIR [rev v2] z] \ \ \ \ \ \ \ \ \ \ \ \ \ \ \ \ \ \
\ \ \ \ \ \ \ \ \ \ \ \ \ MP 8 5

10. UNDIR v1 v3 $|$ UNDIR [rev v2] v3 \ \ \ \ \ \ \ \ \ \ \ \ \ \ \ \ \ \ \
\ \ \ \ \ \ \ \ \ \ \ \ \ \ \ US (v3 z) 9

11. UNDIR [rev v2] v3 \ \ \ \ \ \ \ \ \ \ \ \ \ \ \ \ \ \ \ \ \ \ \ \ \ \ \
\ \ \ \ \ \ \ \ \ \ \ \ \ \ \ \ \ \ \ \ \ \ \ \ \ \ \ \ \ \ \ LDS 10 6

12. (Ay)[UNDIR [rev v2] y

-$>$ (Az)[UNDIR [rev v2] z $|$ UNDIR y z]] \ \ \ \ \ \ \ \ \ \ \ \ \ \ \ \ \
\ \ \ \ US (rev(v2) x) 4

13. UNDIR [rev v2] v3

-$>$ (Az)[UNDIR [rev v2] z $|$ UNDIR v3 z] \ \ \ \ \ \ \ \ \ \ \ \ \ \ \ \ \
\ \ \ \ \ \ \ \ \ \ \ \ \ \ US (v3 y) 12

14. (Az)[UNDIR [rev v2] z $|$ UNDIR v3 z] \ \ \ \ \ \ \ \ \ \ \ \ \ \ \ \ \
\ \ \ \ \ \ \ \ \ \ \ \ \ MP 13 11

15. UNDIR [rev v2] [rev v2] $|$ UNDIR v3 [rev v2] \ \ \ \ \ \ \ \ \ \ US
(rev(v2) z) 14

16. (Ev11)UNDIR v11 v11 $|$ UNDIR v3 [rev v2] \ \ \ \ \ \ \ \ \ \ \ \ \ \ \
\ \ \ \ \ EG 15

17. UNDIR v3 [rev v2] \ \ \ \ \ \ \ \ \ \ \ \ \ \ \ \ \ \ \ \ \ \ \ \ \ \ \
\ \ \ \ \ \ \ \ \ \ \ \ \ \ \ \ \ \ \ \ \ \ \ \ \ \ \ \ \ \ \ LDS 2 16

18. (Ay)[UNDIR v3 [rev y] -$>$ UNDIR y [rev v3]] \ \ \ \ \ \ \ \ \ \ \ \ \ \
\ \ \ \ \ \ US (v3 x) 3

19. UNDIR v3 [rev v2] -$>$ UNDIR v2 [rev v3] \ \ \ \ \ \ \ \ \ \ \ \ \ \ \ \
\ \ \ \ \ \ \ \ \ US (v2 y) 18

20. UNDIR v2 [rev v3] \ \ \ \ \ \ \ \ \ \ \ \ \ \ \ \ \ \ \ \ \ \ \ \ \ \ \
\ \ \ \ \ \ \ \ \ \ \ \ \ \ \ \ \ \ \ \ \ \ \ \ \ \ \ \ \ \ \ \ MP 19 17

21. UNDIR v2 [rev v3] \ \ \ \ \ \ \ \ \ \ \ \ \ \ \ \ \ \ \ \ \ \ \ \ \ \ \
\ \ \ \ \ \ \ \ \ \ \ \ \ \ \ \ \ \ \ \ \ \ \ \ \ \ \ \ \ \ \ \ SAME 20

22. \symbol{126}UNDIR v1 v3 -$>$ UNDIR v2 [rev v3] \ \ \ \ \ \ \ \ \ \ \ \ \
\ \ \ \ \ \ \ \ \ \ \ \ \ \ \ \ \ \ CP 21

23. UNDIR v1 v3 $|$ UNDIR v2 [rev v3] \ \ \ \ \ \ \ \ \ \ \ \ \ \ \ \ \ \ \
\ \ \ \ \ \ \ \ \ \ \ \ \ \ \ \ \ IMP 22

24. UNDIR v1 [rev v2]

-$>$ UNDIR v1 v3 $|$ UNDIR v2 [rev v3] \ \ \ \ \ \ \ \ \ \ \ \ \ \ \ \ \ \ \
\ \ \ \ \ \ \ \ \ \ \ \ \ \ \ \ \ \ \ \ CP 23

25. [\symbol{126}UNDIR v1 [rev v2] $|$ UNDIR v1 v3] $|$ UNDIR v2 [rev v3] \
\ \ IMP 24

26. UNDIR v1 v2

-$>$ [\symbol{126}UNDIR v1 [rev v2] $|$ UNDIR v1 v3] $|$ UNDIR v2 [rev v3] \
\ \ \ \ \ CP 25

27. [[\symbol{126}UNDIR v1 v2 $|$ \symbol{126}UNDIR v1

[rev v2]] $|$ UNDIR v1 v3] $|$ UNDIR v2 [rev v3] \ \ \ \ \ \ \ \ \ \ \ \ \ \
\ \ \ \ \ \ \ \ \ \ \ \ \ \ \ IMP 26

28. (Ax)(Ay)(Az)[[[\symbol{126}UNDIR x y $|$ \symbol{126}UNDIR x [rev y]]

$|$ UNDIR x z] $|$ UNDIR y [rev z]] \ \ \ \ \ \ \ \ \ \ \ \ \ \ \ \ \ \ \ \
\ \ \ \ \ \ \ \ \ \ \ \ \ \ \ \ \ \ \ \ \ \ \ \ \ \ \ \ \ \ \ UG 27

29. [(Ax)\symbol{126}UNDIR x x

\& (Ax)(Ay)[UNDIR x [rev y] -$>$ UNDIR y [rev x]]]

\& (Ax)(Ay)[UNDIR x y -$>$ (Az)[UNDIR x z $|$ UNDIR y z]]

-$>$ (Ax)(Ay)(Az)[[[\symbol{126}UNDIR x y $|$ \symbol{126}UNDIR x [rev y]]

$|$ UNDIR x z] $|$ UNDIR y [rev z]] \ \ \ \ \ \ \ \ \ \ \ \ \ \ \ \ \ \ \ \
\ \ \ \ \ \ \ \ \ \ \ \ \ \ \ \ \ \ \ \ \ \ \ \ \ \ \ \ \ \ \ CP 28

\section{Appendix C. Derive W3 from the axioms I.5, I.6, and ODO.}

Solution

1. [(Ax)\symbol{126}UNDIR x x

\& (Ax)(Ay)[UNDIR x y -$>$ (Az)[UNDIR x z $|$ UNDIR y z]]]

\& (Ax)(Ay)(Az)[\symbol{126}UNDIR x [rev y] \& \symbol{126}UNDIR x z -$>$ 
\symbol{126}UNDIR y [rev z]]

\ \ \ \ \ \ \ \ \ \ \ \ \ \ \ \ \ \ \ \ \ \ \ \ \ \ \ \ \ \ \ \ \ \ \ \ \ \
\ \ \ \ \ \ \ \ \ \ \ \ \ \ \ \ \ \ \ \ \ \ \ \ \ \ \ \ \ \ \ \ \ \ \ \ \ \
\ \ \ \ ASSUMED-PREMISE

2. \symbol{126}(Ex)UNDIR x x \ \ \ \ \ \ \ \ \ \ \ \ \ \ \ \ \ \ \ \ \ \ \ \
\ \ \ \ \ \ \ \ \ \ \ \ \ \ \ \ \ \ \ \ \ \ \ \ \ \ \ \ \ \ \ \ \ \ \ \ \ \
\ SIMP 1

3. (Ax)(Ay)[UNDIR x y-$>$ (Az)[UNDIR x z $|$ UNDIR y z]] \ \ \ \ \ \ \ \ \
SIMP 1

4. (Ax)(Ay)(Az)[\symbol{126}UNDIR x [rev y] \& \symbol{126}UNDIR x z

-$>$ \symbol{126}UNDIR y [rev z]] \ \ \ \ \ \ \ \ \ \ \ \ \ \ \ \ \ \ \ \ \
\ \ \ \ \ \ \ \ \ \ \ \ \ \ \ \ \ \ \ \ \ \ \ \ \ \ \ \ \ \ \ \ \ \ \ \ \ \
\ \ SIMP 1

5. UNDIR v1 [rev v2] \ \ \ \ \ \ \ \ \ \ \ \ \ \ \ \ \ \ \ \ \ \ \ \ \ \ \ \
\ \ \ \ \ \ \ \ \ \ \ \ \ \ \ \ \ \ \ \ \ \ ASSUMED-PREMISE

6. \symbol{126}UNDIR v1 [rev v3] \ \ \ \ \ \ \ \ \ \ \ \ \ \ \ \ \ \ \ \ \ \
\ \ \ \ \ \ \ \ \ \ \ \ \ \ \ \ \ \ \ \ \ \ \ \ \ ASSUMED-PREMISE

7. (Ay)(Az)[\symbol{126}UNDIR v1 [rev y] \& \symbol{126}UNDIR v1 z

-$>$ \symbol{126}UNDIR y [rev z]] \ \ \ \ \ \ \ \ \ \ \ \ \ \ \ \ \ \ \ \ \
\ \ \ \ \ \ \ \ \ \ \ \ \ \ \ \ \ \ \ \ \ \ \ \ \ \ \ \ \ \ US (v1 x) 4

8. (Ay)[UNDIR v1 y

-$>$ (Az)[UNDIR v1 z $|$ UNDIR y z]] \ \ \ \ \ \ \ \ \ \ \ \ \ \ \ \ \ \ \ \
\ \ \ \ \ \ \ \ \ \ \ US (v1 x) 3

9. (Az)[\symbol{126}UNDIR v1 [rev v3] \& \symbol{126}UNDIR v1 z

-$>$ \symbol{126}UNDIR v3 [rev z]] \ \ \ \ \ \ \ \ \ \ \ \ \ \ \ \ \ \ \ \ \
\ \ \ \ \ \ \ \ \ \ \ \ \ \ \ \ \ \ \ \ \ \ \ \ \ \ \ \ US (v3 y) 7

10. UNDIR v1 [rev v2]

-$>$ (Az)[UNDIR v1 z $|$ UNDIR [rev v2] z] \ \ \ \ \ \ \ \ \ \ \ \ \ \ \ \ \
\ \ \ \ \ \ \ US (rev(v2) y) 8

11. (Az)[UNDIR v1 z $|$ UNDIR [rev v2] z] \ \ \ \ \ \ \ \ \ \ \ \ \ \ \ \ \
\ \ \ \ \ MP 10 5

12. UNDIR v1 [rev v3] $|$ UNDIR [rev v2] [rev v3] \ \ \ \ \ \ \ \ \ \ \ US
(rev(v3) z) 11

13. UNDIR [rev v2] [rev v3] \ \ \ \ \ \ \ \ \ \ \ \ \ \ \ \ \ \ \ \ \ \ \ \
\ \ \ \ \ \ \ \ \ \ \ \ \ \ \ \ \ \ LDS 12 6

14. (Ay)[UNDIR [rev v2] y

-$>$ (Az)[UNDIR [rev v2] z $|$ UNDIR y z]] \ \ \ \ \ \ \ \ \ \ \ \ \ \ \ \ \
\ \ \ \ \ \ \ \ US (rev(v2) x) 3

15. UNDIR [rev v2] [rev v3]

-$>$ (Az)[UNDIR [rev v2] z $|$ UNDIR [rev v3] z] \ \ \ \ \ \ \ \ \ \ \ \ \ \
\ \ US (rev(v3) y) 14

16. (Az)[UNDIR [rev v2] z $|$ UNDIR [rev v3] z] \ \ \ \ \ \ \ \ \ \ \ \ \ \
\ MP 15 13

17. UNDIR [rev v2] [rev v2] $|$ UNDIR [rev v3] [rev v2] \ \ \ \ \ US
(rev(v2) z) 16

18. (Ev7)UNDIR v7 v7 $|$ UNDIR [rev v3] [rev v2] \ \ \ \ \ \ \ \ \ \ \ \ \
EG 17

19. UNDIR [rev v3] [rev v2] \ \ \ \ \ \ \ \ \ \ \ \ \ \ \ \ \ \ \ \ \ \ \ \
\ \ \ \ \ \ \ \ \ \ \ \ \ \ \ \ \ \ \ \ \ LDS 2 18

20. \symbol{126}UNDIR v1 [rev v3] \& \symbol{126}UNDIR v1 [rev v3]

-$>$ \symbol{126}UNDIR v3 [rev [rev v3]] \ \ \ \ \ \ \ \ \ \ \ \ \ \ \ \ \ \
\ \ \ \ \ \ \ \ \ \ \ \ \ \ US (rev(v3) z) 9

21. UNDIR v1 [rev v3] $|$ [UNDIR v1 [rev v3]

$|$ \symbol{126}UNDIR v3 [rev [rev v3]]] \ \ \ \ \ \ \ \ \ \ \ \ \ \ \ \ \ \
\ \ \ \ \ \ \ \ \ \ \ \ \ \ \ \ \ \ \ \ \ \ \ \ \ \ \ \ \ IMP 20

22. UNDIR v1 [rev v3] $|$ \symbol{126}UNDIR v3 [rev [rev v3]] \ \ \ \ \ \ \
\ \ \ \ \ \ LDS 21 6

23. \symbol{126}UNDIR v3 [rev [rev v3]] \ \ \ \ \ \ \ \ \ \ \ \ \ \ \ \ \ \
\ \ \ \ \ \ \ \ \ \ \ \ \ \ \ \ \ \ \ \ \ \ \ \ \ \ \ LDS 22 6

24. (Ay)(Az)[\symbol{126}UNDIR v3 [rev y] \& \symbol{126}UNDIR v3 z

-$>$ \symbol{126}UNDIR y [rev z]] \ \ \ \ \ \ \ \ \ \ \ \ \ \ \ \ \ \ \ \ \
\ \ \ \ \ \ \ \ \ \ \ \ \ \ \ \ \ \ \ \ \ \ \ \ \ \ \ \ \ \ \ \ \ \ \ \ US
(v3 x) 4

25. (Ay)[UNDIR v3 y

-$>$ (Az)[UNDIR v3 z $|$ UNDIR y z]] \ \ \ \ \ \ \ \ \ \ \ \ \ \ \ \ \ \ \ \
\ \ \ \ \ \ \ \ \ \ \ \ \ \ \ \ \ US (v3 x) 3

26. (Az)[\symbol{126}UNDIR v3 [rev [rev v3]] \& \symbol{126}UNDIR v3 z

-$>$ \symbol{126}UNDIR [rev v3] [rev z]] \ \ \ \ \ \ \ \ \ \ \ \ \ \ \ \ \ \
\ \ \ \ \ \ \ \ \ \ \ \ \ \ \ \ US (rev(v3) y) 24

27. \symbol{126}UNDIR v3 [rev [rev v3]] \& \symbol{126}UNDIR v3 v2

-$>$ \symbol{126}UNDIR [rev v3] [rev v2] \ \ \ \ \ \ \ \ \ \ \ \ \ \ \ \ \ \
\ \ \ \ \ \ \ \ \ \ \ \ \ \ \ \ \ \ \ \ \ \ \ \ \ \ \ \ \ \ US (v2 z) 26

28. UNDIR v3 [rev [rev v3]] $|$

[UNDIR v3 v2 $|$ \symbol{126}UNDIR [rev v3] [rev v2]] \ \ \ \ \ \ \ \ \ \ \
\ \ \ \ \ \ \ \ \ \ \ \ \ \ \ \ IMP 27

29. UNDIR v3 v2 $|$ \symbol{126}UNDIR [rev v3] [rev v2] \ \ \ \ \ \ \ \ \ \
\ \ \ \ \ \ \ \ \ \ \ \ \ LDS 28 23

30. UNDIR v3 v2 \ \ \ \ \ \ \ \ \ \ \ \ \ \ \ \ \ \ \ \ \ \ \ \ \ \ \ \ \ \
\ \ \ \ \ \ \ \ \ \ \ \ \ \ \ \ \ \ \ \ \ \ \ \ \ \ \ \ \ \ \ \ \ \ RDS 29 19

31. UNDIR v3 v2

-$>$ (Az)[UNDIR v3 z $|$ UNDIR v2 z] \ \ \ \ \ \ \ \ \ \ \ \ \ \ \ \ \ \ \ \
\ \ \ \ \ \ \ \ \ \ \ \ \ \ \ \ \ \ US (v2 y) 25

32. (Az)[UNDIR v3 z $|$ UNDIR v2 z] \ \ \ \ \ \ \ \ \ \ \ \ \ \ \ \ \ \ \ \
\ \ \ \ \ \ \ \ \ \ \ \ \ \ \ \ \ MP 31 30

33. UNDIR v3 v3 $|$ UNDIR v2 v3 \ \ \ \ \ \ \ \ \ \ \ \ \ \ \ \ \ \ \ \ \ \
\ \ \ \ \ \ \ \ \ \ \ \ \ \ \ \ \ \ \ \ US (v3 z) 32

34. (Ev15)UNDIR v15 v15 $|$ UNDIR v2 v3 \ \ \ \ \ \ \ \ \ \ \ \ \ \ \ \ \ \
\ \ \ \ \ \ \ \ \ \ EG 33

35. UNDIR v2 v3 \ \ \ \ \ \ \ \ \ \ \ \ \ \ \ \ \ \ \ \ \ \ \ \ \ \ \ \ \ \
\ \ \ \ \ \ \ \ \ \ \ \ \ \ \ \ \ \ \ \ \ \ \ \ \ \ \ \ \ \ \ \ \ \ \ \ LDS
2 34

36. UNDIR v2 v3 \ \ \ \ \ \ \ \ \ \ \ \ \ \ \ \ \ \ \ \ \ \ \ \ \ \ \ \ \ \
\ \ \ \ \ \ \ \ \ \ \ \ \ \ \ \ \ \ \ \ \ \ \ \ \ \ \ \ \ \ \ \ \ \ \ \ SAME
35

37. \symbol{126}UNDIR v1 [rev v3] -$>$ UNDIR v2 v3 \ \ \ \ \ \ \ \ \ \ \ \ \
\ \ \ \ \ \ \ \ \ \ \ \ \ \ \ \ \ \ CP 36

38. UNDIR v1 [rev v3] $|$ UNDIR v2 v3 \ \ \ \ \ \ \ \ \ \ \ \ \ \ \ \ \ \ \
\ \ \ \ \ \ \ \ \ \ \ \ \ \ \ \ \ IMP 37

39. UNDIR v1 [rev v2]

-$>$ UNDIR v1 [rev v3] $|$ UNDIR v2 v3 \ \ \ \ \ \ \ \ \ \ \ \ \ \ \ \ \ \ \
\ \ \ \ \ \ \ \ \ \ \ \ \ \ \ \ \ \ \ CP 38

40. [\symbol{126}UNDIR v1 [rev v2] $|$ UNDIR v1 [rev v3]] $|$ UNDIR v2 v3 \
\ \ \ \ IMP 39

41. UNDIR v1 v2

-$>$ [\symbol{126}UNDIR v1 [rev v2] $|$ UNDIR v1 [rev v3]] $|$ UNDIR v2 v3 \
\ \ \ \ \ \ CP 40

42. [[\symbol{126}UNDIR v1 v2 $|$ \symbol{126}UNDIR v1 [rev v2]] $|$ UNDIR
v1 [rev v3]]

$|$ UNDIR v2 v3 \ \ \ \ \ \ \ \ \ \ \ \ \ \ \ \ \ \ \ \ \ \ \ \ \ \ \ \ \ \
\ \ \ \ \ \ \ \ \ \ \ \ \ \ \ \ \ \ \ \ \ \ \ \ \ \ \ \ \ \ \ \ \ \ \ \ \ \
\ \ \ \ \ \ \ \ IMP 41

43. (Ax)(Ay)(Az)[[[\symbol{126}UNDIR x y $|$ \symbol{126}UNDIR x [rev y]]

$|$ UNDIR x[rev z]] $|$ UNDIR y z] \ \ \ \ \ \ \ \ \ \ \ \ \ \ \ \ \ \ \ \ \
\ \ \ \ \ \ \ \ \ \ \ \ \ \ \ \ \ \ \ \ \ \ \ \ \ \ \ \ \ \ \ \ UG 42

44. [(Ax)\symbol{126}UNDIR x x

\& (Ax)(Ay)[UNDIR x y -$>$ (Az)[UNDIR x z $|$ UNDIR y z]]]

\& (Ax)(Ay)(Az)[\symbol{126}UNDIR x [rev y] \& \symbol{126}UNDIR x z -$>$ 
\symbol{126}UNDIR y [rev z]]

-$>$ (Ax)(Ay)(Az)[[[\symbol{126}UNDIR x y $|$ \symbol{126}UNDIR x [rev y]] $%
| $ UNDIR x [rev z]]

$|$ UNDIR y z] \ \ \ \ \ \ \ \ \ \ \ \ \ \ \ \ \ \ \ \ \ \ \ \ \ \ \ \ \ \ \
\ \ \ \ \ \ \ \ \ \ \ \ \ \ \ \ \ \ \ \ \ \ \ \ \ \ \ \ \ \ \ \ \ \ \ \ \ \
\ \ \ \ \ \ \ \ \ \ \ \ CP 43

\section{Appendix D. The axiom I.6 implies W4}

Show that

(Ax)(Ay)[UNDIR x y -$>$ (Az)[UNDIR x z $|$ UNDIR y z]]

-$>$ (Ax)(Ay)(Az)[[[\symbol{126}UNDIR x y $|$ \symbol{126}UNDIR x [rev y]]

$|$ UNDIR x [rev z]] $|$ UNDIR y [rev z]]

Solution

1. (Ax)(Ay)[UNDIR x y

-$>$ (Az)[UNDIR x z $|$ UNDIR y z]] \ \ \ \ \ \ \ \ \ \ \ \ \ \ \ \ \ \ \ \
\ \ \ \ \ \ \ \ \ \ \ \ ASSUMED-PREMISE

2. UNDIR v1 v2 \ \ \ \ \ \ \ \ \ \ \ \ \ \ \ \ \ \ \ \ \ \ \ \ \ \ \ \ \ \ \
\ \ \ \ \ \ \ \ \ \ \ \ \ \ \ \ \ \ \ \ \ \ \ \ \ \ ASSUMED-PREMISE

3. \symbol{126}UNDIR v1 [rev v3] \ \ \ \ \ \ \ \ \ \ \ \ \ \ \ \ \ \ \ \ \ \
\ \ \ \ \ \ \ \ \ \ \ \ \ \ \ \ \ \ \ \ \ \ \ \ \ \ ASSUMED-PREMISE

4. (Ay)[UNDIR v1 y -$>$ (Az)[UNDIR v1 z $|$ UNDIR y z]] \ \ \ \ \ \ \ \ \ \
US (v1 x) 1

5. UNDIR v1 v2

-$>$ (Az)[UNDIR v1 z $|$ UNDIR v2 z] \ \ \ \ \ \ \ \ \ \ \ \ \ \ \ \ \ \ \ \
\ \ \ \ \ \ \ \ \ \ \ \ \ \ \ \ \ \ \ \ US (v2 y) 4

6. (Az)[UNDIR v1 z $|$ UNDIR v2 z] \ \ \ \ \ \ \ \ \ \ \ \ \ \ \ \ \ \ \ \ \
\ \ \ \ \ \ \ \ \ \ \ \ \ \ \ \ \ \ \ MP 5 2

7. UNDIR v1 [rev v3] $|$ UNDIR v2 [rev v3] \ \ \ \ \ \ \ \ \ \ \ \ \ \ \ \ \
\ \ US (rev(v3) z) 6

8. UNDIR v2 [rev v3] \ \ \ \ \ \ \ \ \ \ \ \ \ \ \ \ \ \ \ \ \ \ \ \ \ \ \ \
\ \ \ \ \ \ \ \ \ \ \ \ \ \ \ \ \ \ \ \ \ \ \ \ \ \ \ \ \ \ \ \ \ LDS 7 3

9. UNDIR v2 [rev v3] \ \ \ \ \ \ \ \ \ \ \ \ \ \ \ \ \ \ \ \ \ \ \ \ \ \ \ \
\ \ \ \ \ \ \ \ \ \ \ \ \ \ \ \ \ \ \ \ \ \ \ \ \ \ \ \ \ \ \ \ \ SAME 8

10. \symbol{126}UNDIR v1 [rev v3] -$>$ UNDIR v2 [rev v3] \ \ \ \ \ \ \ \ \ \
\ \ \ \ \ \ \ \ \ \ \ \ \ \ CP 9

11. UNDIR v1 [rev v3] $|$ UNDIR v2 [rev v3] \ \ \ \ \ \ \ \ \ \ \ \ \ \ \ \
\ \ \ \ \ \ \ \ \ \ \ \ \ IMP 10

12. UNDIR v1 [rev v2]

-$>$ UNDIR v1 [rev v3] $|$ UNDIR v2 [rev v3] \ \ \ \ \ \ \ \ \ \ \ \ \ \ \ \
\ \ \ \ \ \ \ \ \ \ \ \ \ \ \ \ CP 11

13. [\symbol{126}UNDIR v1 [rev v2] $|$ UNDIR v1 [rev v3]]

$|$ UNDIR v2 [rev v3] \ \ \ \ \ \ \ \ \ \ \ \ \ \ \ \ \ \ \ \ \ \ \ \ \ \ \
\ \ \ \ \ \ \ \ \ \ \ \ \ \ \ \ \ \ \ \ \ \ \ \ \ \ \ \ \ \ \ \ \ \ \ \ \ \
IMP 12

14. UNDIR v1 v2

-$>$ [\symbol{126}UNDIR v1 [rev v2] $|$ UNDIR v1 [rev v3]] $|$ UNDIR v2 [rev
v3]

\ \ \ \ \ \ \ \ \ \ \ \ \ \ \ \ \ \ \ \ \ \ \ \ \ \ \ \ \ \ \ \ \ \ \ \ \ \
\ \ \ \ \ \ \ \ \ \ \ \ \ \ \ \ \ \ \ \ \ \ \ \ \ \ \ \ \ \ \ \ \ \ \ \ \ \
\ \ \ \ \ \ \ \ \ \ \ \ \ \ \ \ \ \ CP 13

15. [[\symbol{126}UNDIR v1 v2 $|$ \symbol{126}UNDIR v1 [rev v2]] $|$ UNDIR
v1 [rev v3]]

$|$ UNDIR v2 [rev v3] \ \ \ \ \ \ \ \ \ \ \ \ \ \ \ \ \ \ \ \ \ \ \ \ \ \ \
\ \ \ \ \ \ \ \ \ \ \ \ \ \ \ \ \ \ \ \ \ \ \ \ \ \ \ \ \ \ \ \ \ \ \ \ \ \
\ IMP 14

16. (Ax)(Ay)(Az)[[[\symbol{126}UNDIR x y $|$ \symbol{126}UNDIR x [rev y]]

$|$ UNDIR x [rev z]] $|$ UNDIR y [rev z]] \ \ \ \ \ \ \ \ \ \ \ \ \ \ \ \ \
\ \ \ \ \ \ \ \ \ \ \ \ \ \ \ \ \ \ \ \ \ \ \ \ UG 15

17. (Ax)(Ay)[UNDIR x y -$>$ (Az)[UNDIR x z $|$ UNDIR y z]]

-$>$ (Ax)(Ay)(Az)[[[\symbol{126}UNDIR x y $|$ \symbol{126}UNDIR x [rev y]]

$|$ UNDIR x [rev z]] $|$ UNDIR y [rev z]] \ \ \ \ \ \ \ \ \ \ \ \ \ \ \ \ \
\ \ \ \ \ \ \ \ \ \ \ \ \ \ \ \ \ \ \ \ \ \ \ \ CP 16

\section{Appendix E \ Derive the axiom I.6 from the axioms I.7, I.8, and ODO}

Solution

1. [(Ax)(Ay)[UNDIR x y $|$ UNDIR x [rev y]] \&

(Ax)(Ay)(Az)[\symbol{126}UNDIR x [rev y] \& \symbol{126}UNDIR x z

-$>$ \symbol{126}UNDIR y [rev z]]] \& (Ax)(Ay)[UNDIR x y \& UNDIR x [rev y]

-$>$ (Az)[UNDIR x z \& UNDIR x [rev z] $|$ UNDIR y z \& UNDIR y [rev z]]]

\ \ \ \ \ \ \ \ \ \ \ \ \ \ \ \ \ \ \ \ \ \ \ \ \ \ \ \ \ \ \ \ \ \ \ \ \ \
\ \ \ \ \ \ \ \ \ \ \ \ \ \ \ \ \ \ \ \ \ \ \ \ \ \ \ \ \ \ \ \ \ \ \ \ \ \
\ \ ASSUMED-PREMISE

2. (Ax)(Ay)[UNDIR x y $|$ UNDIR x [rev y]] \ \ \ \ \ \ \ \ \ \ \ \ \ \ \ \ \
\ \ \ \ \ \ \ \ \ \ \ SIMP 1

3. (Ax)(Ay)(Az)[\symbol{126}UNDIR x [rev y] \& \symbol{126}UNDIR x z -$>$ 
\symbol{126}UNDIR y [rev z]]

\ \ \ \ \ \ \ \ \ \ \ \ \ \ \ \ \ \ \ \ \ \ \ \ \ \ \ \ \ \ \ \ \ \ \ \ \ \
\ \ \ \ \ \ \ \ \ \ \ \ \ \ \ \ \ \ \ \ \ \ \ \ \ \ \ \ \ \ \ \ \ \ \ \ \ \
\ \ \ \ \ \ \ \ \ \ \ \ \ \ \ \ SIMP 1

4. (Ax)(Ay)[UNDIR x y \& UNDIR x [rev y]

-$>$(Az)[UNDIR x z \& UNDIR x [rev z] $|$ UNDIR y z \& UNDIR y [rev z]]]

\ \ \ \ \ \ \ \ \ \ \ \ \ \ \ \ \ \ \ \ \ \ \ \ \ \ \ \ \ \ \ \ \ \ \ \ \ \
\ \ \ \ \ \ \ \ \ \ \ \ \ \ \ \ \ \ \ \ \ \ \ \ \ \ \ \ \ \ \ \ \ \ \ \ \ \
\ \ \ \ \ \ \ \ \ \ \ \ \ \ \ \ SIMP 1

5. UNDIR v1 v2 \ \ \ \ \ \ \ \ \ \ \ \ \ \ \ \ \ \ \ \ \ \ \ \ \ \ \ \ \ \ \
\ \ \ \ \ \ \ \ \ \ \ \ \ \ \ \ \ \ \ \ \ \ \ \ \ ASSUMED-PREMISE

6. \symbol{126}UNDIR v1 v3 \ \ \ \ \ \ \ \ \ \ \ \ \ \ \ \ \ \ \ \ \ \ \ \ \
\ \ \ \ \ \ \ \ \ \ \ \ \ \ \ \ \ \ \ \ \ \ \ \ \ \ \ \ \ ASSUMED-PREMISE

7. (Ay)[UNDIR v1 y \& UNDIR v1 [rev y]

-$>$(Az)[UNDIR v1 z \& UNDIR v1 [rev z] $|$ UNDIR y z \& UNDIR y [rev z]]]

\ \ \ \ \ \ \ \ \ \ \ \ \ \ \ \ \ \ \ \ \ \ \ \ \ \ \ \ \ \ \ \ \ \ \ \ \ \
\ \ \ \ \ \ \ \ \ \ \ \ \ \ \ \ \ \ \ \ \ \ \ \ \ \ \ \ \ \ \ \ \ \ \ \ \ \
\ \ \ \ \ \ US (v1 x) 4

8. UNDIR v1 v2 \& UNDIR v1 [rev v2]

-$>$(Az)[UNDIR v1 z \& UNDIR v1 [rev z]

$|$ UNDIR v2 z \& UNDIR v2 [rev z]] \ \ \ \ \ \ \ \ \ \ \ \ \ \ \ \ \ \ \ \
\ \ \ \ \ \ \ \ US (v2 y) 7

9. \symbol{126}UNDIR v1 v2 $|$ [\symbol{126}UNDIR v1 [rev v2] $|$

(Az)[UNDIR v1 z \& UNDIR v1 [rev z]

$|$ UNDIR v2 z \& UNDIR v2 [rev z]]] \ \ \ \ \ \ \ \ \ \ \ \ \ \ \ \ \ \ \ \
\ \ \ \ \ \ \ \ \ \ IMP 8

10. \symbol{126}UNDIR v1 [rev v2] $|$

(Az)[UNDIR v1 z \& UNDIR v1 [rev z] $|$ UNDIR v2 z

\& UNDIR v2 [rev z]] \ \ \ \ \ \ \ \ \ \ \ \ \ \ \ \ \ \ \ \ \ \ \ \ \ \ \ \
\ \ \ \ \ \ \ \ \ \ \ \ \ \ \ \ \ \ \ \ \ \ \ \ \ LDS 9 5

11. \symbol{126}UNDIR v1 [rev v2] \ \ \ \ \ \ \ \ \ \ \ \ \ \ \ \ \ \ \ \ \
\ \ \ \ \ \ \ \ \ \ \ \ \ \ \ \ \ \ \ \ \ \ \ \ \ \ CASE2 10

12. (Az)[UNDIR v1 z \& UNDIR v1 [rev z] $|$ UNDIR v2 z

\& UNDIR v2 [rev z]] \ \ \ \ \ \ \ \ \ \ \ \ \ \ \ \ \ \ \ \ \ \ \ \ \ \ \ \
\ \ \ \ \ \ \ \ \ \ \ \ \ \ \ \ \ \ \ \ \ \ \ \ \ CASE1 10

13. UNDIR v1 v3 \& UNDIR v1

[rev v3] $|$ UNDIR v2 v3 \& UNDIR v2 [rev v3] \ \ \ \ \ \ \ \ \ \ \ \ \ \ \
\ US (v3 z) 12

14. [UNDIR v1 v3 \& UNDIR v1

[rev v3] $|$ UNDIR v2 v3] \& [UNDIR v1 v3 \& UNDIR v1

[rev v3] $|$ UNDIR v2 [rev v3]] \ \ \ \ \ \ \ \ \ \ \ \ \ \ \ \ \ \ \ \ \ \
DISTRIBUTIVE-LAW 13

15. UNDIR v1 v3 \& UNDIR v1 [rev v3] $|$ UNDIR v2 v3 \ \ \ \ \ \ SIMP 14

16. [UNDIR v1 v3 $|$ UNDIR v2 v3] \& [UNDIR v1

[rev v3] $|$ UNDIR v2 v3] \ \ \ \ \ \ \ \ \ \ \ \ \ \ \ \ \ \ \ \ \ \ \ \ \
\ \ \ \ DISTRIBUTIVE-LAW 15

17. UNDIR v1 v3 $|$ UNDIR v2 v3 \ \ \ \ \ \ \ \ \ \ \ \ \ \ \ \ \ \ \ \ \ \
\ \ \ \ \ \ \ \ \ \ \ \ SIMP 16

18. UNDIR v2 v3 \ \ \ \ \ \ \ \ \ \ \ \ \ \ \ \ \ \ \ \ \ \ \ \ \ \ \ \ \ \
\ \ \ \ \ \ \ \ \ \ \ \ \ \ \ \ \ \ \ \ \ \ \ \ \ \ LDS 17 6

19. (Ay)(Az)[\symbol{126}UNDIR v1 [rev y] \& \symbol{126}UNDIR v1 z

-$>$ \symbol{126}UNDIR y [rev z]] \ \ \ \ \ \ \ \ \ \ \ \ \ \ \ \ \ \ \ \ \
\ \ \ \ \ \ \ \ \ \ \ \ \ \ \ \ \ \ \ \ \ \ \ \ \ \ \ \ \ \ US (v1 x) 3

20. (Az)[\symbol{126}UNDIR v1 [rev v2] \& \symbol{126}UNDIR v1 z

-$>$ \symbol{126}UNDIR v2 [rev z]] \ \ \ \ \ \ \ \ \ \ \ \ \ \ \ \ \ \ \ \ \
\ \ \ \ \ \ \ \ \ \ \ \ \ \ \ \ \ \ \ \ \ \ \ \ \ \ \ \ \ \ \ US (v2 y) 19

21. \symbol{126}UNDIR v1 [rev v2] \& \symbol{126}UNDIR v1 v3

-$>$ \symbol{126}UNDIR v2 [rev v3] \ \ \ \ \ \ \ \ \ \ \ \ \ \ \ \ \ \ \ \ \
\ \ \ \ \ \ \ \ \ \ \ \ \ \ \ \ \ \ \ \ \ \ \ \ \ \ \ \ \ US (v3 z) 20

22. UNDIR v1 [rev v2] $|$

[UNDIR v1 v3 $|$ \symbol{126}UNDIR v2 [rev v3]] \ \ \ \ \ \ \ \ \ \ \ \ \ \
\ \ \ \ \ \ \ \ \ \ \ \ \ IMP 21

23. UNDIR v1 v3 $|$ \symbol{126}UNDIR v2 [rev v3]\ \ \ \ \ \ \ \ \ \ \ \ \ \
\ \ \ \ \ \ \ LDS 22 11

24. \symbol{126}UNDIR v2 [rev v3] \ \ \ \ \ \ \ \ \ \ \ \ \ \ \ \ \ \ \ \ \
\ \ \ \ \ \ \ \ \ \ \ \ \ \ \ \ \ \ \ \ \ \ \ \ \ \ \ LDS 23 6

25. (Ay)(Az)[\symbol{126}UNDIR v2 [rev y] \& \symbol{126}UNDIR v2 z

-$>$ \symbol{126}UNDIR y [rev z]] \ \ \ \ \ \ \ \ \ \ \ \ \ \ \ \ \ \ \ \ \
\ \ \ \ \ \ \ \ \ \ \ \ \ \ \ \ \ \ \ \ \ \ \ \ \ \ \ \ \ \ US (v2 x) 3

26. (Az)[\symbol{126}UNDIR v2 [rev v3] \& \symbol{126}UNDIR v2 z

-$>$ \symbol{126}UNDIR v3 [rev z]] \ \ \ \ \ \ \ \ \ \ \ \ \ \ \ \ \ \ \ \ \
\ \ \ \ \ \ \ \ \ \ \ \ \ \ \ \ \ \ \ \ \ \ \ \ \ \ \ \ \ US (v3 y) 25

27. \symbol{126}UNDIR v2 [rev v3] \& \symbol{126}UNDIR v2 [rev v3]

-$>$ \symbol{126}UNDIR v3 [rev [rev v3]] \ \ \ \ \ \ \ \ \ \ \ \ \ \ \ \ \ \
\ \ \ \ \ \ \ \ \ \ \ \ \ \ \ \ \ \ \ \ \ \ \ US (rev(v3) z) 26

28. UNDIR v2 [rev v3] $|$ [UNDIR v2

[rev v3] $|$ \symbol{126}UNDIR v3 [rev [rev v3]]] \ \ \ \ \ \ \ \ \ \ \ \ \
\ \ \ \ \ \ \ \ \ \ \ \ \ \ \ \ \ IMP 27

29. UNDIR v2 [rev v3] $|$ \symbol{126}UNDIR v3 [rev [rev v3]] \ \ \ \ \ \ \
\ \ LDS 28 24

30. \symbol{126}UNDIR v3 [rev [rev v3]] \ \ \ \ \ \ \ \ \ \ \ \ \ \ \ \ \ \
\ \ \ \ \ \ \ \ \ \ \ \ \ \ \ \ \ \ \ \ \ \ LDS 29 24

31. (Ay)[UNDIR v3 y $|$ UNDIR v3 [rev y]] \ \ \ \ \ \ \ \ \ \ \ \ \ \ \ \ \
\ \ \ \ US (v3 x) 2

32. (Ay)(Az)[\symbol{126}UNDIR v3 [rev y] \& \symbol{126}UNDIR v3 z

-$>$ \symbol{126}UNDIR y [rev z]] \ \ \ \ \ \ \ \ \ \ \ \ \ \ \ \ \ \ \ \ \
\ \ \ \ \ \ \ \ \ \ \ \ \ \ \ \ \ \ \ \ \ \ \ \ \ \ \ \ \ US (v3 x) 3

33. UNDIR v3 [rev [rev v3]] $|$ UNDIR v3 [rev[rev [rev v3]]] \ \ \ 

\ \ \ \ \ \ \ \ \ \ \ \ \ \ \ \ \ \ \ \ \ \ \ \ \ \ \ \ \ \ \ \ \ \ \ \ \ \
\ \ \ \ \ \ \ \ \ \ \ \ \ \ \ \ \ \ \ \ \ \ \ \ \ \ \ \ \ \ \ \ \ US
(rev(rev(v3)) y) 31

34. UNDIR v3 [rev [rev [rev v3]]] \ \ \ \ \ \ \ \ \ \ \ \ \ \ \ \ \ \ \ \ \
\ \ \ \ \ \ \ \ \ \ \ LDS 33 30

35. (Az)[\symbol{126}UNDIR v3 [rev v3] \& \symbol{126}UNDIR v3 z

-$>$ \symbol{126}UNDIR v3 [rev z]] \ \ \ \ \ \ \ \ \ \ \ \ \ \ \ \ \ \ \ \ \
\ \ \ \ \ \ \ \ \ \ \ \ \ \ \ \ \ \ \ \ \ \ \ \ \ \ \ \ \ \ US (v3 y) 32

36. \symbol{126}UNDIR v3 [rev v3] \& \symbol{126}UNDIR v3 [rev

[rev v3]] -$>$ \symbol{126}UNDIR v3 [rev [rev [rev v3]]]

\ \ \ \ \ \ \ \ \ \ \ \ \ \ \ \ \ \ \ \ \ \ \ \ \ \ \ \ \ \ \ \ \ \ \ \ \ \
\ \ \ \ \ \ \ \ \ \ \ \ \ \ \ \ \ \ \ \ \ \ \ \ \ \ \ \ \ \ \ \ \ \ \ \ \ \
\ US (rev(rev(v3)) z) 35

37. \symbol{126}[\symbol{126}UNDIR v3 [rev v3] \& \symbol{126}UNDIR v3 [rev
[rev v3]]] \ \ \ \ MT 36 34

38. UNDIR v3 [rev v3] $|$ UNDIR v3 [rev [rev v3]] \ \ \ \ \ \ \ \ DE.MORGAN
37

39. UNDIR v3 [rev v3] \ \ \ \ \ \ \ \ \ \ \ \ \ \ \ \ \ \ \ \ \ \ \ \ \ \ \
\ \ \ \ \ \ \ \ \ \ \ \ \ \ \ \ \ \ \ \ RDS 38 30

40. \symbol{126}UNDIR v2 [rev v3] \& \symbol{126}UNDIR v2 v3

-$>$ \symbol{126}UNDIR v3 [rev v3] \ \ \ \ \ \ \ \ \ \ \ \ \ \ \ \ \ \ \ \ \
\ \ \ \ \ \ \ \ \ \ \ \ \ \ \ \ \ \ \ \ \ \ \ US (v3 z) 26

41. UNDIR v2 [rev v3] $|$

[UNDIR v2 v3 $|$ \symbol{126}UNDIR v3 [rev v3]] \ \ \ \ \ \ \ \ \ \ \ \ \ \
\ \ \ \ \ \ \ \ IMP 40

42. UNDIR v2 v3 $|$ \symbol{126}UNDIR v3 [rev v3] \ \ \ \ \ \ \ \ \ \ \ \ \
\ \ \ \ \ \ LDS 41 24

43. UNDIR v2 v3 \ \ \ \ \ \ \ \ \ \ \ \ \ \ \ \ \ \ \ \ \ \ \ \ \ \ \ \ \ \
\ \ \ \ \ \ \ \ \ \ \ \ \ \ \ \ \ \ \ \ \ \ \ RDS 42 39

44. UNDIR v2 v3 \ \ \ \ \ \ \ \ \ \ \ \ \ \ \ \ \ \ \ \ \ \ \ \ \ \ \ \ \ \
\ \ \ \ \ \ \ \ \ \ \ \ \ \ \ \ \ \ \ \ \ \ \ SAME 18

45. UNDIR v2 v3 \ \ \ \ \ \ \ \ \ \ \ \ \ \ \ \ \ \ \ \ \ \ \ \ \ \ \ \ \ \
\ \ \ \ \ \ \ \ \ \ \ \ \ \ \ \ \ \ \ \ \ \ \ SAME 43

46. UNDIR v2 v3 \ \ \ \ \ \ \ \ \ \ \ \ \ \ \ \ \ \ \ \ \ \ \ \ \ \ \ \ \ \
\ \ \ \ \ \ \ \ \ \ \ \ \ \ \ \ \ \ \ \ \ \ \ CASES 10 45 44

47. \symbol{126}UNDIR v1 v3 -$>$ UNDIR v2 v3 \ \ \ \ \ \ \ \ \ \ \ \ \ \ \ \
\ \ \ \ \ \ \ \ \ \ CP 46

48. UNDIR v1 v3 $|$ UNDIR v2 v3 \ \ \ \ \ \ \ \ \ \ \ \ \ \ \ \ \ \ \ \ \ \
\ \ \ \ \ \ \ \ IMP 47

49. (Az)[UNDIR v1 z $|$ UNDIR v2 z] \ \ \ \ \ \ \ \ \ \ \ \ \ \ \ \ \ \ \ \
\ \ \ \ \ UG 48

50. UNDIR v1 v2

-$>$ (Az)[UNDIR v1 z $|$ UNDIR v2 z] \ \ \ \ \ \ \ \ \ \ \ \ \ \ \ \ \ \ \ \
\ \ \ \ \ \ CP 49

51. (Ax)(Ay)[UNDIR x y -$>$ (Az)[UNDIR x z $|$ UNDIR y z]] \ \ \ \ UG 50

52. [(Ax)(Ay)[UNDIR x y $|$ UNDIR x [rev y]]

\&(Ax)(Ay)(Az)[\symbol{126}UNDIR x [rev y] \& \symbol{126}UNDIR x z -$>$ 
\symbol{126}UNDIR y [rev z]]]

\& (Ax)(Ay)[UNDIR x y \& UNDIR x [rev y]

-$>$(Az)[UNDIR x z \& UNDIR x [rev z] $|$ UNDIR y z \& UNDIR y [rev z]]]

-$>$ (Ax)(Ay)[UNDIR x y -$>$ (Az)[UNDIR x z $|$ UNDIR y z]] \ \ \ \ \ \ \ \
\ CP 51


\begin{thebibliography}{99}
\bibitem{Heyting} Heyting, A.: \ Axioms for intuitionistic plane affine
geometry In L. Henkin et al., eds.. The Axiomatic Method, North-Holland,
Amsterdam, 160-173 (1959).

\bibitem{jvp-95} von Plato, J.: The axioms of constructive geometry. Annals
of Pure and Applied Logic, vol. 76, 169-200 (1995).

\bibitem{Dalen} van Dalen, D.: `Outside' as a primitive notion in
constructive projective geometry. Geometriae Dedicata 60, 107- 111 (1996).

\bibitem{jvp-98} von Plato, J.: A constructive theory of ordered affine
geometry. Indag. Mathem., N.S. 9 (4), 549-562 (1998).

\bibitem{Wos} McCune, W., Wos, L.: Application of automated deduction to the
search for single axioms for exponent groups. Proc. of LPAR-92, St.
Petersburg, Russia, July 15--20, 1992, pp. 131--136.

\bibitem{McCune-93} Padmanabhan, R., McCune, W.: Single identities for
ternary Boolean algebra. Comput. Math. Appl., 29 (2), 13-16 (1993).

\bibitem{McCune-96} McCune, W., Padmanabhan, R.: Single identities for
lattice theory and weakly associative lattices. Algebra Universalis, 36 (4),
436-449(1996).

\bibitem{McCune} McCune, W.: Solution of the Robbins Problem. J. of
automated reasoning 19(3), 263--276 (1997).

\bibitem{dli-97} Li, D.: Unification Algorithms for Eliminating and
Introducing Quantifiers in Natural Deduction Automated Theorem Proving. J.
of Automated Reasoning 18(1), 105-134(1997).

\bibitem{dli-tableau} Li, D.: Natural Deduction Prover and Experiments,
Tableaux 97, Nancy, France, pp153-157, May 13-16, (1997).

\bibitem{dli-halt} Li, D.: A Mechanical Proof of the Halting Problem in
Natural Deduction Style, AAR Newsletter No.23, June 1993.

\bibitem{dli-00} Li, D., Jia, P., Li, X.: Simplifying von Plato's
axiomatization of constructive apartness geometry. Annals of pure and
applied logic, vol. 102, No1-2, (2000).

\bibitem{dli-03} Li, D.: Using Prover ANDP to simplify orthogonality. Annals
of pure and applied logic, 124, 49-70 (2003).

\bibitem{dli-04} Li, D.: The equality axioms are not independent. ACM SIGACT
NEWS, Vol. 35, Issue 3, 98-101 (2004).
\end{thebibliography}
\end{document}